\documentclass[a4paper,11pt,leqno]{amsart}
\usepackage{amssymb,hyperref}
\usepackage[latin1]{inputenc}
\usepackage[english]{babel}
\usepackage{graphicx}
\usepackage{color}

\usepackage[T1]{fontenc}
\usepackage[latin1]{inputenc}
\usepackage[all]{xy}
\usepackage{stmaryrd}
\usepackage{epsfig,psfig}
\usepackage{float}

\usepackage{psfrag}
\usepackage{wrapfig}
\usepackage{amsmath}
\usepackage{amssymb,latexsym}
\usepackage[mathscr]{eucal}
\usepackage{url}
\usepackage{fancybox}
\swapnumbers

\numberwithin{equation}{section}

\DeclareMathAlphabet{\mathrmsl}{OT1}{cmr}{m}{sl}

\newcommand{\oper}[3][n]{\newcommand{#2}{\mathop{\mathrm{#3}}%
\ifx n#1\nolimits\else\limits\fi} }
\newcommand{\rsoper}[3][n]{\newcommand{#2}{\mathop{\mathrmsl{#3}}%
\ifx n#1\nolimits\else\limits\fi} }

\newcounter{mnotecount}[section]
\renewcommand{\themnotecount}{\thesection.\arabic{mnotecount}}
\newcommand{\mnote}[1]
{\protect{\stepcounter{mnotecount}}$^{\mbox{\footnotesize  $
      \bullet$\themnotecount}}$ \marginpar{\raggedright\tiny\em
    $\!\!\!\!\!\!\,\bullet$\themnotecount: #1} }

\newcommand{\eq}[1]{\ref{#1}}

\renewcommand{\qed}{\hfill $\blacksquare$ \medskip \\}

\newcommand{\question}{\noindent {\bf Question:\ }}

\theoremstyle{main}
\newtheorem{thm} {\bf  Theorem} [section]
\newtheorem{coro} [thm] {\bf  Corollary}
\newtheorem{lem} [thm] {\bf  Lemma}
\newtheorem{prop} [thm] {\bf  Proposition}

\newtheorem{d\'ef}[thm]{\bf Definition}
\newtheorem{Rq}[thm]{\bf  Remark }

\DeclareFontFamily{OT1}{rsfs}{}
\DeclareFontShape{OT1}{rsfs}{m}{n}{ <-7> rsfs5 <7-10> rsfs7 <10-> rsfs10}{}
\DeclareMathAlphabet{\mycal}{OT1}{rsfs}{m}{n}

\newcommand{\dVol}{\operatorname{dVol}}
\newcommand{\Vol}{\operatorname{Vol}}
\newcommand{\Scal}{\operatorname{Scal}}
\newcommand{\Int}{\operatorname{Int}}
\newcommand{\Ric}{\operatorname{Ric}}
\newcommand{\tr}{ {\rm tr}}

\newcommand{\hess}{\operatorname{Hess}}
\newcommand{\Ker}{\operatorname{Ker}}
\renewcommand{\d}{{\rm d}}

\newcommand{\Span}{\operatorname{Span}}

\newcommand{\so}{\operatorname{so}}

\newcommand{\R}{{\mathbb R}}

\newcommand{\Sn}{{\mathbb S^{n}}}

\title[KID on totally umbilical \& compact hypersurfaces]{Killing initial data on totally umbilical \& compact hypersurfaces}
\author{Daniel MAERTEN}
\date{\today}
\address{Université Paris 6, Institut de Mathématiques de Jussieu, UMR CNRS 7586, 175 rue du Chevaleret, 75013 Paris, France} 
\keywords{Scalar curvature, harmonic curvature, Einstein metrics, Schwarzschild--de Sitter metrics, constraints application, Lafontaine equation}
\email{{\color{red}daniel.maerten@yahoo.fr}}
\thanks{2000 \textit{Mathematics Subject Classification.} 53C24, 53C25, 53C80, 83C05, 83C15}
\thanks{The author benefits from the ANR grant \textit{Géométrie des variétés d'Einstein non compactes ou singulières}, n° ANR--06--BLAN--0154--02}

\begin{document}

\begin{abstract}
	In this note, we give a geometric characterization of the compact and totally umbilical hypersurfaces that carry a non trivial locally static Killing Initial Data (KID). More precisely, such compact hypersurfaces $(M^n,g, cg)$ endowed with a Riemannian metric $g$ and a second fundamental form $cg$ (where $c\in C^{\infty}(M)$ a priori) have constant mean curvature and are isometric to one of the following manifolds:
\begin{itemize}
	\item[(i)] $\Sn$ the standard sphere,
	\item[(ii)] a finite quotient  of a warped product $(\mathbb S^1 \times Y, \d t^2 + h^2 (t)g_0)$, where $(Y^{n-1},g_0)$ is Einstein with positive scalar curvature.
\end{itemize}
In particular, they have harmonic curvature and strictly positive constant scalar curvature.
\end{abstract}

\maketitle

\section{Introduction}

A good starting point for this article is the study of the scalar curvature application
$$ u : \mycal{M} \longrightarrow C^\infty (M) , \ g \longmapsto \Scal^g \ ,$$
where $\mycal M$ denotes the cone of the Riemannian metrics on $M$. If we denote by $U_g$ the differential of $u$ at a metric $g\in\mycal M$, then for any $h\in \Gamma(S^2 T^*M)$ we have (see \cite{B} for details)
$$U_g(h)= \Delta (\tr_g (h)) + \delta (\delta h) - \left\langle \Ric^g, h\right\rangle \ ,$$
where $\Delta$ is the positive Laplacian with respect to $g$, $\left\langle \cdot,\cdot\right\rangle$ is the inner product extended to tensors of any type and $\delta$ is the $g$--divergence operator defined as
$$ \delta S (X_1,\cdots,X_p)= -\sum^{n}_{i=1}\nabla_{e_i}S(e_i,X_1, \cdots,X_p) \ ,$$
for any $(p+1)$--tensor $S$ $(p\ge 0)$ and for any local orthonormal $g$--basis $(e_i)^{n}_{i=1}$. The formal adjoint of $U_g$ is given by 
$$\forall f \in C^\infty \quad U^*_g(f)= \hess^g f -f \Ric^g + (\Delta f) g \ .$$
From a geometric point of view, the fact that $\Ker U^*_g\neq \left\{0\right\}$ is a necessary condition for the scalar curvature application $u$ not to be a submersion at $g$, and so a necessary condition for the level set $u^{-1}(\Scal^g)\subset\mycal M$ not to be a submanifold (in the sense of infinite dimensional submanifolds, see once again \cite{B}). A natural issue is to study the compact and connected manifolds $(M^n,g)$ that admit a non trivial function solution to the equation
$$ \hess^g f -f \Ric^g + (\Delta f) g =0 \qquad (*) \ ,$$ 
which has been extensively studied (notably by Lafontaine, that is why we will call $(*)$ the Lafontaine equation). As a matter of fact, Fischer and Marsden \cite{FM} believed that the standard sphere $\Sn$ was the unique compact and connected Riemannian manifold (except the Ricci flat compact manifolds) that have non trivial solutions to the Lafontaine equation \cite{FM}. Their conjecture is false, since Lafontaine listed in \cite{L} all the compact and conformally flat manifolds admitting solutions  to $(*)$, in particular some non Ricci flat, non spherical warped product metrics appear in that list. However, without the conformally flatness assumption, we do not have to our disposal an exhaustive list of manifolds that have non trivial solutions of $(*)$. This question still remains open nowadays. Another interesting geometric interpretation of $(*)$ that is particularly relevant from General Relativity, deals with the classification of static solutions of the Einstein equations. More precisely, if $\gamma=-f^2\d t^2 +g$ is a static Lorentzian warped product metric on $\R\times M$, then the metric $\gamma$ is Einstein (which is equivalent to be a solution of the Einstein vacuum equations) if and only if $f$ is a solution of $(*)$. Thus, determining the solutions of $(*)$ is strictly equivalent to classify the static solutions of Einstein vacuum equations. This fact is closely related to the Lorentzian metrics obtained thanks to a Killing development (see below for details).\\
An important feature is that $(*)$ reduces (up to a normalization) to Obata's equation $\hess^g f =-f  g $ when the metric $g$ is assumed to be Einstein, and therefore it is a generalization of Obata's equation. It is well known that the existence of a non trivial function solution of Obata's equation characterizes the geometry of the standard sphere $\Sn$ \cite{O}.\\
When we fix a metric $g_0\in \mycal M$ and restrict $u$ to the space $\mycal V$ of metrics that have the same Riemannian volume form than $g_0$, namely $\mycal V= \left\{g\in\mycal M/ \dVol_g= \dVol_{g_0 }\right\}$, then if $u$ is not a submersion at $g$ then $\Ker \left(U^*_g\right)_0\neq \left\{0\right\}$  with $\left(U^*_g\right)_0$ the traceless part of $U^*_g$ namely
$$\forall f \in C^\infty \quad \left(U^*_g(f)\right)_0= \hess^g f -f \Ric^g_0 + \frac{\Delta f}{n} g \ ,$$
where $\Ric_0^g$ denotes the traceless part of the Ricci curvature. In \cite{L}, Lafontaine exhibited a large family of metrics that admit non trivial solutions to $\hess^g f -f \Ric^g_0 + \frac{\Delta f}{n} g=0$.\\
One can also be interested in the (positive) constants that could be the scalar curvature of a metric $g$ among the metrics of a certain fixed volume, let us say $\Vol(\Sn)$. This is equivalent to restrict $u$ to the space $\mycal C =  \left\{g\in\mycal M/ \d\Scal_g= 0 \textrm{ and } \Vol(M,g)=\Vol(\Sn)\right\}$. In that case, if $u$ is not a submersion at $g$ then there exists a non trivial solution $f$ of $U^*_g(f)=\Ric^g_0$. All these results are summarized in the following statement.
\begin{thm} Let $M^n$  be a compact and connected manifold $(n\ge 3)$.
\begin{enumerate}
	\item (\textbf{Obata} \cite{O}) If u is not a submersion at an \textbf{Einstein} metric $g\in\mycal M$ then $(M^n,g)$ is isometric to a standard sphere and $\Ker U^*_g=\Span\{x_1,x_2,\cdots, x_{n+1}\}$ with $(x_i)^{n+1}_{i=1}$ the standard coordinates on $\Sn$.
	\item (\textbf{Lafontaine} \cite{L}) If u is not a submersion at a \textbf{conformally flat} metric $g\in\mycal M$ then $(M^n,g)$ is isometric to one of the following manifolds:
\begin{enumerate}
	\item [1)] A flat compact manifold and $\Ker U^*_g= \R$,
	\item [2)] The standard sphere $\Sn$ and $\Ker U^*_g=\Span\{x_1,x_2,\cdots, x_{n+1}\}$.
	\item [3)] A finite quotient of $(\mathbb S^1 \times\mathbb S^{n-1}, \d t^2 + g_{\mathbb S^{n-1}})$ and $\dim\Ker U^*_g=2$.
	\item [4)] A finite quotient of a warped product $(\mathbb S^1 \times\mathbb S^{n-1}, \d t^2 +h^2(t) g_{\mathbb S^{n-1}})$ and ${\Ker U^*_g=\R h'}$.
\end{enumerate}
	\item (\textbf{Lafontaine} \cite{L}) If $Y^{n-1}$ is  a compact and connected manifold that carries an Einstein metric of positive scalar curvature, then there exists an infinite dimensional set in $\mycal V(\mathbb S^1\times Y)=\mycal V$ where $u_{|\mycal V}$ is not a submersion. 
	\item (\textbf{Lafontaine} \cite{L}) If $u_{|\mycal C}$ is not a submersion at a \textbf{conformally flat} metric $g\in \mycal C$ then $(M^n,g)$ is isometric to the standard sphere $\mathbb S^n$.
	\item (\textbf{Bessières--Lafontaine--Rozoy} \cite{BLR}) If $n=3$ and $u_{|\mycal C}$ is not a submersion at a metric $g\in \mycal C$ then $(M^3,g)$ is isometric to the standard sphere $\mathbb S^3$.
\end{enumerate}
\end{thm}

A natural generalization of these questions is obtained by considering the natural extension (due to the General Relativity context) of the scalar curvature application $u$: the constraints application $\Phi$. More precisely, if we consider $(M^n,g,k)\subset (N^{n+1},\gamma)$ a Riemannian submanifold in a Lorentzian manifold, endowed with the induced metric $g$ and the second fundamental form $k$, the constraints application is given by
$$\Phi:(g,k)\longmapsto \left(
\begin{array}{c}
	\Scal^g + (\tr_g k)^2 - \left|k\right|^2_g \\
	-2(\delta k + \d \tr_g k)
\end{array}
\right)= \left(
\begin{array}{c}
	\Phi_1 (g,k) \\
	\Phi_2 (g,k)
\end{array}
\right) \in C^\infty(M)\times\Gamma(T^* M) ,$$
where $(g,k)\in \mycal{M}\times \Gamma(S^2T^* M)$. When the Lorentzian manifold $(N^{n+1},\gamma)$ is a solution of the Einstein equations (it is not necessarily the case for our problem) that is to say when the metric $\gamma$ satisfies $$\Ric^\gamma -\frac{1}{2}\Scal^\gamma \gamma=T\ ,$$
where $T$ is the stress--energy tensor, then the Hamiltonian constraint $\Phi_1(g,k)$ corresponds to $T(\partial_t,\partial_t)$ and the moment constraint $\Phi_2(g,k)$ corresponds to (up to a multiplicative constant) $T(\partial_t,\cdot)_{|TM}$ (here $\partial_t$ denotes a unit normal to $M\hookrightarrow N$). The other pieces  of the tensor $T$ are usually called the Einstein evolution equations and are characterized by the occurrences of the $\partial_t$ derivatives of the metric $\gamma$ (on the contrary to the constraints equations). The constraints equations are also the traced Gauss and Codazzi equations of the embedding $M\hookrightarrow N$.\\
We denote by $L_{(g,k)}$ ($=L$ in short) the differential of $\Phi$ at $(g,k)$, and by $L^*$ its formal adjoint. Analogously to the scalar curvature application $u$, $\Phi$ is not a submersion at some couple $(g,k)$ if and only if $\Ker L^*_{(g,k)}\neq\left\{0\right\}$ . The formal adjoint of $L_{(g,k)}$ is given by the following rather complicated coupled differential operator 
$$\left\{
\begin{array}{lll}
	L_1^* (f,\alpha) & = & E(f,\alpha) - \left(\tr_g E(f,\alpha)\right)g - \frac{1}{2}\Big(\left\langle L_2^* (f,\alpha), k\right\rangle  + \left\langle \alpha, \Phi_2 (g,k)\right\rangle + 2f\Phi_1 (g,k)  \Big)g    \\
	L_2^* (f,\alpha) & = &  -2(\delta^* \alpha +fk) +2 \tr_g(\delta^* \alpha +fk) g
\end{array}\right. \ , $$
where $(f,\alpha)\in C^{\infty}(M) \times \Gamma(T^*M)$ and
$$ E(f,\alpha):= \hess^g f -f(\Ric^g +2 (\tr_g k)k -2k\circ k) + \mycal{L}_\alpha k + (\delta \alpha)k \ . $$
Here $\mycal L$ stands for the Lie derivative, $k\circ k$ means $(k\circ k)_{ij}= k_{ir}k_{sj}g^{rs}$ and finally $\delta^{*}\alpha$ is the symmetric part of the covariant derivative $\nabla\alpha$ i.e. $\delta^{*}\alpha(X,Y)=\frac{1}{2}\big(\nabla_X \alpha(Y)+ \nabla_Y \alpha(X)\big)$. 
The condition $\Ker L^*_{(g,k)}\neq\left\{0\right\}$ is equivalent to the existence of a non trivial couple ${(f,\alpha)\in  C^{\infty}(M) \times \Gamma(T^*M)} $ such that
$$\left\{
\begin{array}{l}
	\hess^g f + \mycal{L}_\alpha k = f(\Ric^g + (\tr_g k)k -2k\circ k) - \frac{1}{2(n-1)}\Big(  \left\langle \alpha, \Phi_2 (g,k)\right\rangle + 2f\Phi_1 (g,k)  \Big)g \\
	\mycal{L}_\alpha g + 2fk = 0
\end{array}\right. \ .$$
In this context we can address the following issue concerning the existence of Killing Initial Data (KID, which are by definition \cite{M} the non trivial elements of $\Ker L^*_{(g,k)}$ or equivalently the non trivial solutions of the differential system above).\\

\question Does the existence of a non trivial KID (i.e. an element in $\Ker L^*_{(g,k)}$) characterize the geometry of the Riemannian hypersurface $(M^n,g,k)$ ?\\

This problem is very difficult in general, however Beig, Chru\'sciel and Schoen proved in \cite{BChS} that KIDs are non generic for a large class of slices $(M^n,g,k)$. More explicitly, they proved that the family $  \left\{(g,k)\in \mycal M \times \Gamma(S^2 T^*M)/ \  \Ker L^*_{(g,k)}=\left\{0\right\}\right\}$, was an open and dense set for a $C^{k,\beta}\times C^{k,\beta}$ type topology. Unfortunately, they were not able to give the geometry of $(M^n,g,k)$ under the \textit{exceptional} condition that $\Ker L^*_{(g,k)}$ is non zero. The aim of the present paper is to study a particular situation that generalizes the Lafontaine equation, namely the KID equations on a compact and totally umbilical hypersurface. In other words, we consider $(M^n,g,k)$ such that
\begin{itemize}
\item[1)]$ \exists c \in C^{\infty}(M), c\not\equiv 0,  \quad k=cg $, 
\item[2)] $ \exists (f,\alpha)\not\equiv (0,0)  \quad L^*_{(g,k)}(f,\alpha)=0 $ .
\end{itemize}
It is clear that the standard (constant mean curvature and totally umbilical) spherical slices in de Sitter space--time satisfy all these conditions. Analogously to the Fischer and Marsden conjecture, one could think that the conditions 1) and 2) characterize the spherical slices in de Sitter. This conjecture is false once again because of a result owed to Lafontaine. To see that we first need to reformulate our problem. In virtue of of the umbilicity condition 1), the moment constraints becomes $\Phi_2 (g, cg)= -2(n-1)\d c$. The Hamiltonian one forces the scalar curvature of $g$ to satisfy: $ \Phi_1(g,cg) =\Scal^g + n(n-1)c^2 $. The second KID equation is exactly  $\mycal{L}_\alpha g + 2cfg= 0 $, which means that $\alpha$ is a Killing conformal 1--form (it is non isometric as soon as $f\not \equiv 0$). Finally, the first KID equation is equivalent to 
\begin{eqnarray*}
	\hess^g f &=& f(\Ric^g +nc^2 g -2 c^2 g ) - \mycal{L}_\alpha (cg) \\
	 && - \frac{1}{2(n-1)}\Big(  -2(n-1)\left\langle \alpha, \d c\right\rangle + 2f (\Scal^g + n(n-1)c^2)  \Big)g \\
	\hess^g f &=& f(\Ric^g +nc^2 g -2 c^2 g ) -  c\mycal{L}_\alpha g - \nabla_{\alpha}(cg) + \left\langle \alpha, \d c\right\rangle - f\left(\frac{\Scal^g}{n-1}\right)g -nfc^2 g  \\
	\hess^g f &=& f \Ric^g - f\left(\frac{\Scal^g}{n-1}\right)g \  ,
\end{eqnarray*}
which is equivalent to $U_g^* (f)= \hess^g f -f \Ric^g + (\Delta f) g =0 $ where we have used the relation obtained  thanks to the traced equation. As a consequence, the conditions 1) and 2) are equivalent to the system
$$(\Sigma)\quad \left\{
\begin{array}{l}
	 \mycal{L}_\alpha g + 2cfg= 0 \\
	 U_g^* (f)=0
\end{array}\right. $$
The conjecture \textit{à la Fischer and Marsden} now reads as: the only metric $g$ (except the Ricci flat) that has non trivial solutions to $(\Sigma)$ is the standard metric on the sphere $\Sn$. This is clearly false since any compact manifold $(M^n,g)$ having a non zero Killing field $\alpha$ provides a non trivial solution of $(\Sigma)$ of the form $(0,\alpha)$ and whatever the mean curvature $c$ is. Consequently, from now on we define as a \textit{non trivial KID} a solution $(f,\alpha)$ such that $f\not\equiv  0$. Even under this new definition of non trivial KID, the conjecture  \textit{à la Fischer and Marsden} is false in virtue of the following result (which is a corollary of classification results of Derdzinski about compact manifolds with harmonic curvature \cite{D1,D2}).
\begin{prop}[Lafontaine \cite{L}]\label{jaclaf}
Let $(M^n,g)$ be a compact manifold $(n\ge 3)$ with harmonic curvature which carries a closed Killing conformal and non isometric 1-form, then  ${\Ker U^*_g \not = \left\{0 \right\}}$ and $(M^n,g)$ is isometric to one of the following manifolds:
\begin{itemize}
	\item[(i)] $\Sn$ the standard sphere,\\
	\item[(ii)] a finite quotient  of a warped product $(\mathbb S^1 \times Y, \d t^2 + h^2 (t)g_0)$, where $(Y^{n-1},g_0)$ is Einstein with positive scalar curvature.
\end{itemize}
\end{prop}
Indeed, under the assumptions of Proposition~\eq{jaclaf}, the harmonic curvature and the second Bianchi identity compel the scalar curvature $\Scal^g$ to be a (positive) constant. We denote by $\alpha$ the closed Killing conformal and non isometric 1--form i.e. $\d \alpha =0$ and ${\mycal{L}_\alpha g + \varphi g= 0}$ (this equation was studied by Tashiro \cite{T}) for a certain function $\varphi\not \equiv 0$. It is proved in \cite{L} that $ U_g^* (\varphi)=0$, and thereby the couple $(\alpha, f:= \frac{\varphi}{2c} )$ for any constant $c\in \R^*$, is a solution of $(\Sigma)$ as soon as $(M^n,g)$ is isometric to one of the manifolds listed in Proposition~\eq{jaclaf}. Moreover, the slice $(M^n,g,cg)$ is constant mean curvature and satisfies the vacuum constraints equations with the positive cosmological constant $\Lambda=\frac{1}{2}\left(\Scal^g + n(n-1)c^2\right)>0$.\\
The first main result of this article is to consider a system $(\Sigma_1)$ slightly stronger than $(\Sigma)$, since we moreover demand the Killing conformal form to be closed, and give the geometry of the manifolds that carry non trivial solutions of this system.
\begin{thm}\label{KID}
Let $(M^n, g)$ be a compact and connected manifold $(n\ge 3)$ with $C^3$ metric, which has a non trivial solution $(f,\alpha)$ of
$$(\Sigma_1)\quad \left\{
\begin{array}{l}
	 \nabla\alpha+ cfg= 0 \\
	 U_g^* (f)=0
\end{array}\right. \ ,$$
for a certain $c\in C^{\infty}(M), c\not\equiv 0$. Then c is a non--zero constant and $(M^n, g)$ is isometric to one of the following manifolds:
\begin{itemize}
	\item[(i)] $\Sn$ the standard sphere,
	\item[(ii)] a finite quotient  of a warped product $(\mathbb S^1 \times Y, \d t^2 + h^2 (t)g_0)$, where $(Y^{n-1},g_0)$ is Einstein with positive scalar curvature.
\end{itemize}
From the KID point of view, the vacuum slice $(M^n, g, cg)$ has positive constant scalar curvature and constant mean curvature c given by $\Scal^g +n(n-1)c^2= \Phi_1(g,cg)=const.$.
\end{thm}
The original definition of the KID is owed to Beig and Chru\'sciel \cite{BCh} and is the following. Consider a Killing field $X$ on the ambient Lorentzian manifold $(N,\gamma)$ and decompose it in terms of a lapse function $f$ and a shift 1--form (or vector) $\alpha$. Then the Killing equation of $X$ restricted to the Riemannian hypersurface $(M^n,g,k)$ gives birth to a system of equations satisfied by the couple $(f,\alpha)$ which is by definition a \textit{Killing Initial Data}. Now when the constraints are vacuum (possibly with a cosmological constant) then being a KID in the original sense of \cite{BCh} is equivalent to belong to $\Ker L^*_{(g,k)}$. The existence of a KID allows us to make a Killing development of $(M^n,g,k)$ i.e. to construct a Lorentzian manifold with topology $\R\times M=\tilde{N}$ which is endowed with the metric
$$ \tilde{\gamma}= (\left|\alpha\right|^2-f^2)\d t^2 +  2 \d t\odot \alpha + g \ ,$$ 
where the functions $f, \alpha_i , g_{ij}$ are trivially extended along the product  $\R\times M$, and $\d t\odot \alpha$ denotes the symmetric part of $\d t\otimes \alpha$. By construction, the metric $\tilde{\gamma}$ is stationary and the spatial slices $\left\{t=const.\right\}$ are isometric to our initial manifold $(M^n,g,k)$. The important fact is that the Lorentzian manifold $(\tilde{N},\tilde{\gamma})$ is a solution of the vacuum Einstein equations (possibly with a cosmological constant) that carries a Killing field $\partial_t$. In fact, $\partial_t$ extends the vector field $(f\nu+ \alpha^i\partial_i) \in \Gamma(T\tilde{N}_{|M})$ (where $\nu$ denotes a unit normal to the slice $\left\{t=0\right\}\cong (M^n,g,k)$) along the whole Lorentzian manifold $(\tilde{N},\tilde{\gamma})$. As already said, the space--time $(\tilde{N},\tilde{\gamma})$ is stationary by construction, but when we assume that the spatial part $\alpha$ is closed, it forces the Killing development to be locally static, in the sense of the local integrability of the orthogonal distribution of $\partial_t$. In particular, by making a Killing development of the metrics that are classified in the result above, we obtain a class of static solutions of the vacuum Einstein equations with a positive cosmological constant: the Schwarzschild--de Sitter metrics. These metrics had been studied in \cite{Bir} with $\Lambda<0$ (the Schwarzschild--anti de Sitter metrics), but all the computations of Birmingham can be carried out with $\Lambda>0$. As a conclusion, we can claim that the Schwarzschild--de Sitter solutions are characterized by the existence of non trivial solutions of $(\Sigma_1)$, and we will call the metrics listed in the theorem above the spatial Schwarzschild--de Sitter metrics ($\Sn$ corresponds to the spatial de Sitter metrics).\\

The second main result concerns the restriction of the constraints application to the set $\mycal V \times \Gamma(S^2 T^*M)$ which gives rise to the problem of finding some $(f,\alpha)$ such that
$$\left\{
\begin{array}{l}
	 \mycal{L}_\alpha g + 2cfg= 0 \\
	 \hess^g f -f \Ric^g_0 + \frac{\Delta f}{n} g =0
\end{array}\right. \Longleftrightarrow  
\left\{
\begin{array}{l}
	 \mycal{L}_\alpha g + 2cfg= 0 \\
	 U_g^* (f)=\frac{1}{n}\Big((n-1) \Delta f -f \Scal^g\Big) g
\end{array}\right. .$$
It is pointless to prove that the function $c$ is a constant as in the situation of the system $(\Sigma_1)$, since there exist some examples (the standard sphere $\Sn$ is one of them) of manifolds $(M^n,g)$ that have non trivial solutions $(f,\alpha)$ with a non constant functions $c$. We will discuss this point later on. However, we can classify the solutions of this system when we assume that $c$ is a non--zero constant and that the Killing conformal form $\alpha$ is closed. We surprisingly obtain the same family of manifolds than for the system $(\Sigma_1)$.
\begin{thm}\label{KID2}
Let $(M^n, g)$ be a compact and connected manifold $(n\ge 3)$ with $C^3$ metric, which has a non trivial solution $(f,\alpha)$ of
$$(\Sigma_2)\quad \left\{
\begin{array}{l}
	 \nabla\alpha+ cfg= 0 \\
	 U_g^* (f)=\frac{1}{n}\Big((n-1) \Delta f -f \Scal^g\Big) g
\end{array}\right. \ ,$$
for a certain $c\in\R^*$. Then $(M^n, g)$ is isometric to one of the following manifolds:
\begin{itemize}
	\item[(i)] $\Sn$ the standard sphere,
	\item[(ii)] a finite quotient  of a warped product $(\mathbb S^1 \times Y, \d t^2 + h^2 (t)g_0)$, where $(Y^{n-1},g_0)$ is Einstein with positive scalar curvature.
\end{itemize}
From the KID point of view, the slice $(M^n, g, cg)$ is a spatial Schwarzschild--de Sitter metric i.e. it is vacuum and has positive constant scalar curvature and constant mean curvature c given by $\Scal^g +n(n-1)c^2= \Phi_1(g,cg)=const.$.
\end{thm}

The third main result concerns the restriction of the constraints application to the set $\mycal C \times \Gamma(S^2 T^*M)$ which gives rise to the problem of finding some $(f,\alpha)$ such that
$$\left\{
\begin{array}{l}
	 \mycal{L}_\alpha g + 2cfg= 0 \\
	 U_g^* (f)=\Ric^g_0
\end{array}\right. \ ,$$
which must be related to the result in dimension 3 of Bessières--Lafontaine--Rozoy evocated earlier. Here again, we study this system assuming the closed character of the Killing conformal form $\alpha$.
\begin{thm}\label{BLR}
Let $(M^n, g)$ be a compact and connected manifold $(n\ge 2)$ with $C^3$ metric and constant scalar curvature. Suppose there exists a non trivial solution $(f,\alpha)$ of
$$(\Sigma_3)\quad \left\{
\begin{array}{l}
	 \nabla\alpha+ cfg= 0 \\
	 U_g^* (f)=\Ric^g_0
\end{array}\right. \ ,$$
$c\in C^{\infty}(M), c\not\equiv 0$. Then c is a non--zero constant and $(M^n,g)$ is isometric to the standard sphere $\Sn$. From the KID point of view, the slice $(M^n, g)$ is a spatial de Sitter metric.
\end{thm}

This article is organized as follows: in Section~\eq{preliminaires} we give some technical preliminary results that will be used in Section~\eq{preuves} in order to prove the main Theorems. Section~\eq{preuves} also contains specific results in small dimension $n=2$ or $3$. Finally, Section~\eq{4} is devoted to the study of another interesting system of equations.


\section{Preliminary results}\label{preliminaires}
 
In this article, $(M^n,g)$ is a compact and connected $n$--dimensional manifold with Levi--Civita connection $\nabla$. A 1--form $\alpha$ is said to be Killing conformal and non isometric if it satisfies the equation $\mycal{L}_\alpha g=\psi g $ for a certain function $\psi\not\equiv 0$. Such a form $\alpha$ is not closed in general, but if it is, then $\alpha$ and $\d \psi$ are closely related.
\begin{prop}
Let $(M^n, g)$ be a compact and connected manifold $(n\ge 2)$. Suppose there exist a non trivial couple $(\psi,\alpha)$ such that $\mycal{L}_\alpha g=\psi g$. If the Killing conformal $\alpha$ is closed then $\alpha \wedge \d \psi=0$.
\end{prop}
\proof The covariant derivative $\nabla\alpha\in \Gamma(\otimes^2 T^*M)$ is the sum of a symmetric part and a skew symmetric part since $\nabla\alpha= \frac{1}{2} \big(\d \alpha + \mycal{L}_\alpha g\big)$, where
\begin{equation*}
	   \left\{
\begin{array}{lll}
	\d \alpha(X,Y) &=&   \nabla_X \alpha(Y) - \nabla_Y \alpha(X) = \d^\nabla \alpha (X,Y) \\
	\mycal{L}_\alpha g (X,Y)&=&  \nabla_X \alpha(Y) + \nabla_Y \alpha(X) = 2\delta^\ast \alpha (X,Y) 
\end{array}\right. \ ,
\end{equation*}
and $\d^{\nabla}S(X,Y,X_1,\cdots,X_p):= \nabla_X S(Y,X_1,\cdots,X_p)- \nabla_Y S(X,X_1,\cdots,X_p)$ for any ${(p+1)}$--tensor $S$. The condition $\d^2 \alpha =0$ implies
\begin{equation*}
	\nabla_Z \d \alpha (X,Y) = - \d^\nabla \d \alpha (X,Y,Z) \ .
\end{equation*}
We apply the Ricci identity to $\alpha$
\begin{eqnarray*}
	R(X,Y,Z,\alpha) &=&  \nabla^2_{X,Y} \alpha(Z) -  \nabla^2_{Y,X} \alpha(Z) \\
	&=& \d^\nabla (\nabla\alpha) (X,Y,Z) \\
	&=& \frac{1}{2}\d^\nabla (\d \alpha) (X,Y,Z)  + \d^\nabla (\delta^\ast \alpha) (X,Y,Z) \\
	&=& - \frac{1}{2}\nabla_Z \d \alpha (X,Y) +  \d^\nabla(\psi g)(X,Y,Z) \\
	&=& - \frac{1}{2}\nabla_Z \d \alpha (X,Y) +  \d \psi \wedge g (X,Y,Z) \ ,
\end{eqnarray*}
where $(\omega \wedge S) (X,Y,Z):= \omega(X)S(Y,Z)- \omega(Y)S(X,Z)$, for every $\omega\in \Gamma(T^*M)$ and every $S\in\Gamma(S^2T^*M)$. By plugging $Z=\alpha$ in this relation and assuming that $\alpha$ is closed, we immediately get $ \d \psi \wedge \alpha=0$.\qed

We can wonder if we can do without the closed assumption in the Lafontaine Proposition~\eq{jaclaf}. The proposition above seems to show that we indeed cannot. The following lemma enumerates equations that will be useful in the sequel.
\begin{lem}\label{lem1}
Let $(M^n, g)$ be a compact and connected manifold $(n\ge 2)$ that admits a closed Killing conformal and non isometric 1--form $\alpha$, i.e. $ \nabla \alpha= \psi g $  for a certain $\psi\in C^{\infty}(M),\ {\psi \not\equiv 0}$. \\
Then the following equations hold:
\begin{eqnarray}
	\label{t1} R(X,Y,Z,\alpha)&=&  \d \psi \wedge g(X,Y,Z) \\
	\label{t2} \d \psi \wedge \alpha &=&0 \\
	\label{t3} \Ric^g (Y,\alpha)&=& -(n-1) \d \psi(Y)\\
	\label{t4} \nabla_X \Ric^g (Y,\alpha)&=& -\Big( (n-1)U_g^* (\psi)  +n\psi\Ric^g -(n-1)(\Delta \psi) g\Big)(X,Y) \\
	\label{t5} \nabla_\alpha \Ric^g(X,Y) &=& -\Big((n-2)U_g^* (\psi) + n\psi \Ric^g -(n-1)(\Delta \psi) g\Big)(X,Y)\\
	\label{t6} \d ^\nabla \Ric^g (\alpha,X,Y)&=& U_g^*(\psi)(X,Y) \\ 
	\label{t7} \d ^\nabla \Ric^g (X,Y,\alpha)&=& 0 \ .
\end{eqnarray}
\end{lem}
\proof Formula~(\eq{t1}-\eq{t2})  were proved in the previous proposition. Equation~(\eq{t3}) is obtained by taking the trace in (\eq{t1}). Equation~(\eq{t4}) is the covariant derivative of (\eq{t3}), namely 
$$\nabla_X \Ric^g (Y,\alpha)= -(n-1) \hess^g \psi(X,Y) -\psi\Ric^g (X,Y) \ ,$$
where we have expressed $\hess^g \psi$ in terms of $U_g^*(\psi)$. Equation~(\eq{t5}) is a corollary of the first variation formula of the Ricci curvature (cf. formula (d) of Theorem 1.174 in \cite{Besse} with $\mycal{L}_\alpha g = 2\psi g$) i.e.
$$ \nabla_\alpha \Ric^g (X,Y)= -(n-2)\hess^g \psi(X,Y) - 2\psi \Ric^g (X,Y) + (\Delta \psi)g(X,Y) \ ,$$
where we have expressed again $\hess^g \psi$ in terms of $U_g^*(\psi)$. Finally, Equation~(\eq{t6}) is the result of the difference between (\eq{t5}) and (\eq{t4}), and Equation~(\eq{t7}) is the skew symmetric part of (\eq{t4}).\qed

It is important to notice that we neither need to fix the value of $U_g^* (\psi)$ nor to suppose that the scalar curvature is a constant to have the formulas of Lemma~\eq{lem1}. The only assumption is to have a closed Killing conformal and non isometric 1--form on the manifold. The manifold $(M^n,g)$ is said to have harmonic curvature if $\d ^\nabla \Ric^g=0$. In view of (\eq{t7}), it seems natural to ask if the existence of a closed Killing conformal and non isometric 1--form implies that $g$ has harmonic curvature. It is not true in general, but we can nonetheless deduce some interesting information on $\d ^\nabla \Ric^g$.

\begin{lem}\label{lem2}
Let $(M^n, g)$ be a compact and connected manifold $(n\ge 2)$ that admits a closed Killing conformal and non isometric 1--form $\alpha$, i.e. $ \nabla \alpha= \psi g $ for a certain  $\psi\in C^{\infty}(M),\ \psi \not\equiv 0$. \\
Then the following equation holds:
\begin{equation}\label{eq}
		\psi\d ^\nabla U_g^* (\psi) = \d \psi \wedge U_g^*(\psi) -\psi^2 \d ^\nabla\Ric^g +\Big( \psi\d (\Delta \psi) - (\Delta \psi) \d \psi\Big)\wedge g \ .
\end{equation} 
\end{lem}

\proof To prove Equation~(\eq{eq}) we take the covariant derivative of (\eq{t1}) 
\begin{eqnarray}\label{e1}
	\nabla_X R(Y,Z,T,\alpha) &=&  \Big\{ \hess^g \psi(X,Y)g(Z,T) - \hess^g f(X,Z)g(Y,T) \Big\} \\
\nonumber	&&  -\psi R(Y,Z,T,X) \ .
\end{eqnarray}
We work on the expression $U_g^* (\psi)$, by computing $\d ^\nabla U_g^* (\psi)$ i.e.
\begin{equation*}
	\d ^\nabla U_g^* (\psi)=  \d ^\nabla \hess^g \psi + \d (\Delta \psi) \wedge g -\psi \d ^\nabla\Ric^g - \d \psi \wedge \Ric^g \ .
\end{equation*}
The Ricci identity applied to $\d \psi$ gives $\d ^\nabla \hess^g \psi(X,Y,Z)= R(X,Y,Z,\nabla \psi)$, and therefore
\begin{equation}\label{e2}
	 \d ^\nabla U_g^* (\psi)= R(\cdot,\cdot,\cdot, \nabla \psi)  + \d (\Delta \psi) \wedge g -\psi \d ^\nabla\Ric^g - \d \psi \wedge \Ric^g \ .
\end{equation}
Besides if we set $T=\alpha$ in (\eq{e1}) and use $ \nabla \alpha  =\psi g $, we find (thanks to the symmetry of the curvature tensor)
\begin{eqnarray*}
		\nabla_X R(Y,Z,\alpha,\alpha)&=& X\cdot \underbrace{R(Y,Z,\alpha,\alpha)}_{=0} - \underbrace{R(\nabla_X Y,Z,\alpha,\alpha)}_{=0} - \underbrace{R(Y,\nabla_X Z,\alpha,\alpha)}_{=0}\\
		&& -\psi \big( \underbrace{R(Y,Z,X,\alpha) +  R(Y,Z,\alpha, X)}_{=0}\big) \ .
\end{eqnarray*}
Thanks to the relation $\alpha \wedge \d \psi=0$, we can use once again (\eq{e1}) and the previous formula to get
\begin{eqnarray*}
	0=\nabla_X R(Y,Z,\nabla \psi,\alpha)&=&   \Big\{ \hess^g \psi(X,Y)g(Z,\nabla \psi) - \hess^g \psi(X,Z)g(Y,\nabla \psi) \Big\} \\
	&&  -\psi R(Y,Z,\nabla \psi,X)\\
	&=& - \d \psi \wedge \hess^g \psi (Y,Z,X) +\psi R(Y,Z,X,\nabla \psi) \ ,
\end{eqnarray*}
that we write in short $\psi R(\cdot,\cdot,\cdot,\nabla \psi)= \d \psi \wedge \hess^g \psi $. Finally, by multiplying (\eq{e2}) by $\psi$ and using our curvature formula it comes out
\begin{eqnarray*}
	\psi\d ^\nabla U_g^* (\psi) &=&  \d \psi \wedge \hess^g \psi + \psi \d (\Delta \psi) \wedge g - \psi \d \psi \wedge \Ric^g -\psi^2 \d ^\nabla\Ric^g \\
	&=& \d \psi \wedge U_g^*(\psi) -\psi^2 \d ^\nabla\Ric^g +\Big( \psi\d (\Delta \psi) - (\Delta \psi) \d \psi\Big)\wedge g \ . 
\end{eqnarray*}
\qed

\begin{Rq} Unfortunately, Formula~(\eq{eq}) is more complicated than we could have expected.  Nevertheless, when $\psi$ is an eigenfunction for the Riemannian Laplacian $\Delta$ i.e. $\Delta \psi= \lambda \psi$ for a constant $\lambda\ge 0$ then $\psi\d (\Delta \psi) - (\Delta \psi) \d \psi= 0 $, which clearly simplifies (\eq{eq}).
\end{Rq}
 
The following result gives additional information when the scalar curvature is supposed to be a constant.

\begin{thm}\label{main}
Let $(M^n, g)$ be a compact and connected manifold $(n\ge 2)$ with constant scalar curvature that admits a closed Killing conformal and non isometric 1--form $\alpha$, i.e. ${\nabla \alpha= \psi g} $ for a certain $\psi\in C^{\infty}(M),\ \psi \not\equiv 0$. \\
Then we have:
\begin{enumerate}
	\item $\Delta \psi = \frac{\Scal^g}{n-1}\psi$ , 
	\item $\Scal^g >0$ ,
	\item $\psi\d ^\nabla U_g^* (\psi) = \d \psi \wedge U_g^*(\psi) -\psi^2 \d ^\nabla\Ric^g$.
\end{enumerate}
\end{thm}
\proof (i) The trace of Equation~(\eq{t5}) leads to 
$$ \d \Scal^g(\alpha)=2(n-1) \Delta \psi -2\psi\Scal^g=0 \ ,$$
since the scalar curvature is a constant and (i) immediately follows.\\
(ii) Thanks to (i) we know that $\frac{\Scal^g}{n-1}$ is an eigenvalue for the Riemannian Laplacian on a compact and connected manifold and consequently $\Scal^g\ge 0$. In particular, $\Scal^g=0$ if and only if the eigenfunction $\psi$ is a constant. But it is in fact impossible. By contradiction, suppose $\d\psi=0=\frac{1}{nc}\d \delta \alpha$ on $M$, so $\d \delta \alpha=0$. By integrating by parts we find 
$$0=\int_{M}\left\langle\d \delta \alpha,\alpha \right\rangle=\int_{M}\left|\delta\alpha\right|^{2} \ , $$
and thereby $\delta\alpha=0$ on $M$ which is impossible since $\alpha$ is non isometric. We deduce that $\psi$ cannot be a constant and thus $\Scal^g>0$.\\
(iii) This is a simplification of Equation~(\eq{eq}) thanks to (i) and the previous remark.\qed

We finish this section by recalling a classical result of Bourguignon which is omnipresent throughout this article.

\begin{thm}[Bourguignon \cite{B}]\label{Bourg}
	Let $(M^n, g)$ be a compact and connected manifold ${(n\ge 2)}$. Suppose that there exists some $\psi\in \Ker U^*_g \setminus \{0\}$, then $\Scal^g$ is a nonnegative constant and the level set $\psi^{-1}(0)$ is either empty ($\psi$ is a constant) or composed with embedded hypersurfaces of M.
\end{thm}
\proof Firstly, taking the trace of $ U_g^* (\psi)=0$ gives $\Delta\psi=\frac{\Scal^g}{n-1}\psi$, so that it becomes
$\hess^g \psi -\psi\Ric^g + \frac{\Scal^g}{n-1}\psi g=0$. Secondly, take $x\in M$ and $t\mapsto \sigma(t)$ a geodesic curve such that $\sigma(0)=x$. Consider the function $F:= \psi\circ\sigma$ which has to satisfy the order 2 ordinary differential equation
$$ F''(t) -F(t)\left\{\Ric^g(\sigma'(t),\sigma'(t))- \frac{\Scal^g}{n-1}g(\sigma'(t),\sigma'(t))\right\}=  0  \ ,$$
with the initial conditions $F(0)=\psi(x), \  F'(0)=\d_x \psi(\sigma'(0))$. Now, if one supposes that $x\in\psi^{-1}(0)$ is critical for $\psi$, then the initial conditions vanish and $F$ has to be identically zero. Since we can cover a dense part of $M$ with geodesics starting at $x$, $\psi$ should also vanish on the whole $M$, which is not permitted. We conclude that if $\psi$ is not constant,  $\psi^{-1}(0)$ contains no critical point and so is composed with embedded hypersurfaces of $M$. Besides $\psi$ is a non--zero constant if and only if $g$ is Ricci flat, and in that case $\psi^{-1}(0)=\emptyset$. As regards the scalar curvature, the trick is to compute the divergence of $U^*_g(\psi)=0$
$$0=\delta U_g^* (\psi)=\delta \hess^g \psi +\Ric^g(\nabla \psi)- \psi \delta\Ric^g - \d(\Delta \psi)= \frac{1}{2}\psi \d \Scal^g \ ,$$
where we have used $\delta \Ric^g= -\frac{1}{2}\d \Scal^g$ (cf. 3.135 i) in \cite{GHL}),  $\delta \hess^g f=\d (\Delta f) -\Ric^g(\nabla f)$ (cf. the proof of 4.14 in \cite{GHL}). In any case  $\psi^{-1}(\R^*)$ is dense in $M$, thereby $\d \Scal^g =0$ on the whole manifold. But $\frac{\Scal^g}{n-1}$ is an eigenvalue of the Riemannian Laplacian  on a compact manifold, which implies $\Scal^g\ge 0$. \qed


\section{Proof of the Theorems}\label{preuves}

\subsection{The system $(\Sigma_1)$} We first prove a general result stating that the existence of a non trivial KID on a compact and totally umbilical hypersurface, has to be constant mean curvature in any dimension greater than 2.

\begin{thm}\label{cmc}
	Let  $(M^n, g)$ be a compact and connected manifold $(n\ge2)$ which has a solution $(f,\alpha), f\not \equiv 0$, of
$$ \left\{
\begin{array}{l}
	 \mycal{L}_\alpha g+ 2cfg= 0 \\
	 U_g^* (f)=0
\end{array}\right. \ ,$$
for a certain $c\in C^{\infty}(M), c\not\equiv 0$. Then the scalar curvature $\Scal^g$ is a positive constant and $c$ is a non zero constant.
\end{thm}
\proof Thanks to Theorem~\eq{Bourg}, we know that $\Scal^g$ is a nonnegative constant and that the function $f$ is an eigenfunction since $\Delta f=\frac{\Scal^g}{n-1}f$. The first variation Formula~(\eq{t5}) of the Ricci curvature applies since $\mycal{L}_\alpha g=- 2cfg$:
$$ \nabla_\alpha \Ric^g (X,Y)= -(n-2)\hess^g (cf)(X,Y) - 2cf \Ric^g (X,Y) + (\Delta (cf))g(X,Y) \ ,$$
whose trace is zero since $\Scal^g$ is constant. We derive $\Delta (cf)=\frac{\Scal^g}{n-1}cf$. Suppose now that $\Scal^g=0$ then $fc$ should be constant and so the Killing conformal form $\alpha$ should be isometric, which is in contradiction with $c\not\equiv0, f\not \equiv0$. We have then proved that $\Scal^g>0$, and also that the function $cf$ is an eigenfunction which lies in the same eigenspace than $f$. Now consider a nodal domain of $f$ that we denote by $\Omega$. Using Theorem~\eq{Bourg}, we conclude that $\partial\Omega$ is a regular hypersurface of $M$, and $f$ is a solution of the usual Dirichlet problem on $\Omega$ since by construction we have 
$$\left\{
\begin{array}{cl}
	\Delta f=\frac{\Scal^g}{n-1}f& \textrm{on } \Omega  \\
	f=0 & \textrm{on } \partial\Omega
\end{array}\right. \ ,$$
with positive eigenvalue $\frac{\Scal^g}{n-1}>0$. But, since $\Omega$ is a nodal domain, $f$ cannot vanish and so $\frac{\Scal^g}{n-1}>0$ has to be the first Dirichlet eigenvalue of the domain $(\Omega,g)$ which is always simple. Now, the function $(cf)$ also satisfies
$$\left\{
\begin{array}{cl}
	\Delta (cf)=\frac{\Scal^g}{n-1}cf& \textrm{on } \Omega  \\
	cf=0 & \textrm{on } \partial\Omega
\end{array}\right. \ ,$$
and consequently, there exists a non zero constant $K\in\R^*$ (otherwise $f$ should be zero on an open set which is excluded by Theorem~\eq{Bourg}) such that $Kcf=f$. But $f$ never vanishes on $\Omega$, so that $c=K^{-1}$ has to be a non zero constant on $\Omega$. This  argument is true on any nodal domain of $f$ whose union is a dense open set in $M$, which means that $c$ has to be a non zero constant on the whole $M$.\qed

As a consequence we can completely answer the question of characterizing the totally umbilical KID in dimension 2.
\begin{coro}
	Let  $(M^2, g)$ be a compact and connected Riemannian surface which has a solution $(f,\alpha), f\not\equiv 0$, of
$$ \left\{
\begin{array}{l}
	 \mycal{L}_\alpha g+ 2cfg= 0 \\
	 U_g^* (f)=0
\end{array}\right. \ ,$$
for a certain $c\in C^{\infty}(M), c\not\equiv 0$. Then $c$ is a non zero constant and $(M^2,g)$ is isometric to the standard sphere $\mathbb S^2$. In other words, any totally umbilical and compact hypersurface having a non trivial KID is isometric to a spatial (spherical and constant mean curvature) de Sitter metric.
\end{coro}
\proof From Theorem~\eq{cmc} we know that $c$ is a non zero constant and $\Scal^g$ a positive constant. Gauss--Bonnet formula claims that $M^2$ is homeomorphic to $\mathbb S^2$ and Theorem 3.83 in \cite{GHL} says that $M^2$ is in fact isometric to a standard 2--dimensional sphere (since the scalar curvature and the sectional curvature are equivalent for surfaces).\qed

For manifolds of dimension greater than 3, we need to make an additional assumption (the Killing conformal 1--form $\alpha$ is supposed to be closed) in order to obtain geometric information. We then prove Theorem~\eq{KID} that classifies the compact and connected manifolds $(M^n, g)$ that carry non trivial solutions of $(\Sigma_1)$.

\begin{thm}\label{KID}
Let $(M^n, g)$ be a compact and connected manifold $(n\ge3)$ with $C^3$ metric, which has a solution $(f,\alpha), f\not\equiv0$, of
$$(\Sigma_1)\quad \left\{
\begin{array}{l}
	 \nabla\alpha+ cfg= 0 \\
	 U_g^* (f)=0
\end{array}\right. \ ,$$
for a certain $c\in C^{\infty}(M), c\not\equiv 0$. Then c is a non--zero constant and $(M^n, g)$ is isometric to one of the following manifolds:
\begin{itemize}
	\item[(i)] $\Sn$ the standard sphere,
	\item[(ii)] a finite quotient  of a warped product $(\mathbb S^1 \times Y, \d t^2 + h^2 (t)g_0)$, where $(Y^{n-1},g_0)$ is Einstein with positive scalar curvature.
\end{itemize}
From the KID point of view, the vacuum slice $(M^n, g, cg)$ has positive constant scalar curvature and constant mean curvature c given by $\Scal^g +n(n-1)c^2= \Phi_1(g,cg)=const.$.
\end{thm}
\proof  We can use (iii) in Theorem~\eq{main} with $\psi=cf$ and $c$ a non zero constant. We obtain
$$0=c^2 f\d ^\nabla U_g^* (f) = c^2 \d f \wedge U_g^*(f) -c^2 f^2 \d ^\nabla\Ric^g=  -c^2 f^2 \d ^\nabla\Ric^g \ ,$$
since $U_g^* (f)=0$. It comes out $\d ^\nabla\Ric^g =0 $ on the dense open set $f^{-1}(\R^*)$, and by a continuation argument ($\d ^\nabla\Ric^g \in C^0$ since the metric $g$ is regular enough i.e. $C^3$), we get $\d ^\nabla\Ric^g =0 $ on the whole $M$. We conclude by applying Proposition~\eq{jaclaf} which gives the announced classification.\\
From the KID point of view, $(M^n, g, k)$ is totally umbilical with extrinsic curvature $k=cg$ by construction, and the constant $c\in \R^*$ is the mean curvature of the slice. It is clear that it satisfies the vacuum constraint equations with the positive cosmological constant $\Lambda=\frac{1}{2}\left(\Scal^g + n(n-1)c^2\right)>0$.\qed

\begin{Rq} The function h in the warped product metrics of (ii) in Theorem~\eq{KID} has to satisfy an order 2 differential equation which is given by the relation between  $\Scal^g$ and $\Scal^{g_0}$ the scalar curvature of $(Y,g_0)$. See \cite{L} for further details.
\end{Rq}

In the KID context, Theorems~\eq{KID}  gives rise to the following natural issue.\\
\question Let $(M^n,g)$ be a compact and connected manifold $(n\ge3)$ and ${c\in C^\infty (M), c\not\equiv 0}$. Does the existence of non trivial solutions (i.e. $f\not\equiv0$) of 
$$\left\{
\begin{array}{l}
	 \mycal{L}_\alpha g + 2cfg= 0 \\
	 U_g^* (f)=0
\end{array}\right. ,$$
characterize the spatial Schwarzschild--de Sitter metrics? In general, it is not clear since we a priori lose the important geometric curvature identity (\eq{t1}).

\subsection{The system $(\Sigma_2)$}
As already said in the introduction, we can focus on the restriction of the constraints application to $\mycal V \times \Gamma(S^2 T^*M)$ where  $\mycal V= \left\{g\in\mycal M/ \dVol_g= \dVol_{g_0 }\right\} $. In this context, the relevant operator is $\left(U^*_g\right)_0$ so that we consider the problem of finding some $(f,\alpha)$ such that
$$\left\{
\begin{array}{l}
	 \mycal{L}_\alpha g + 2cfg= 0 \\
	 \hess^g f -f \Ric^g_0 + \frac{\Delta f}{n} g =0
\end{array}\right. \ ,$$
where $c\in C^{\infty}(M), c\not\equiv0$. This problem is clearly equivalent to 
$$\left\{
\begin{array}{l}
	 \mycal{L}_\alpha g + 2cfg= 0 \\
	 U_g^* (f)=\frac{1}{n}\Big((n-1) \Delta f -f \Scal^g\Big) g
\end{array}\right. .$$
The main difference with the system $(\Sigma)$ is that there is no hope to prove that the function $c$ (remind that $c$ is the mean curvature of the hypersurface $(M^n,g,k=cg)$) is a constant. Indeed, let $(M^n,g)$ be an Einstein manifold with nonnegative scalar curvature that has a non isometric Killing conformal 1--form $\alpha$ (such manifolds do exist, the standard sphere $\mathbb S^n$ is one of them). Then any couple $(f,\alpha)$ with $f$ any non zero constant, has to satisfy
$$	\mycal{L}_\alpha g =-2cfg  \quad \textrm{and } \left(U^*_g (f)\right)_0=0  \ ,$$
where we have defined the function $c$ as $c:= \frac{\delta\alpha}{fn}\not\equiv 0$ since $\alpha$ is non isometric. Expecting the non isometric Killing conformal 1--form to be closed does not even guarantee that $c$ should be a constant. You just have to consider on the standard sphere $\mathbb S^n$, the non trivial couples $(f,f\d c)$ where $f$ is any non zero constant and $c\in \Ker U^*_g=\Span\left\{x_1,x_2,\cdots, x_n\right\}$, with $(x_i)^{n+1}_{i=1}$ the standard coordinates on $\mathbb S^n$. Such couples satisfy 
$$\nabla\alpha=-cfg \quad \textrm{and } \left(U^*_g(f)\right)_0=0 \ .$$
Now if we give up the closed character of $\alpha$, then we can consider some non trivial couples  $(f,f\d c + \beta)$ where $f$ is any non zero constant, $c\in \Ker U^*_g=\Span\left\{x_1,x_2,\cdots, x_n\right\}$ and $\beta$ any Killing form  on $\Sn$. Then such couples has to verify 
$$\mycal{L}_\alpha g + 2cfg= 0 \quad \textrm{and } \left(U^*_g(f)\right)_0=0 \ .$$
This space of solutions is parametrized by $\R^*\times \R^{n+1}\times \so(n+1)$ with the Lie algebra $\so(n+1)\cong \R^{\frac{n(n+1)}{2}}$. It is the reason why we restrict our study to the case where $c$ is a non zero constant. Then we can classify the solutions of this system when we assume in addition that the Killing conformal form $\alpha$ is closed, as it is claimed in the following result.

\begin{thm}\label{KID2}
Let $(M^n, g)$ be a compact and connected manifold $(n\ge3)$ with $C^3$ metric, which has a non trivial solution $(f,\alpha), f\not\equiv 0$, of
$$(\Sigma_2)\quad \left\{
\begin{array}{l}
	 \nabla\alpha+ cfg= 0 \\
	 U_g^* (f)=\frac{1}{n}\Big((n-1) \Delta f -f \Scal^g\Big) g
\end{array}\right. \ ,$$
for a certain constant $c \neq 0$. Then $(M^n, g)$ is isometric to one of the following manifolds:
\begin{itemize}
	\item[(i)] $\Sn$ the standard sphere,
	\item[(ii)] a finite quotient  of a warped product $(\mathbb S^1 \times Y, \d t^2 + h^2 (t)g_0)$, where $(Y^{n-1},g_0)$ is Einstein with positive scalar curvature.
\end{itemize}
From the KID point of view, the vacuum slice $(M^n, g, cg)$ has positive constant scalar curvature and constant mean curvature c given by $\Scal^g +n(n-1)c^2= \Phi_1(g,cg)=const.$.
\end{thm}

\proof Our goal is to show that  $U_g^* (f)=0$ and then use Theorem~\eq{KID}. Equation~(\eq{t6}) reads as 
$$\d ^\nabla \Ric^g (\alpha,X,Y)=-c U_g^*(f)(X,Y)=-\frac{c}{n}\Big((n-1) \Delta f -f \Scal^g\Big) g \ . $$
Let $x\in M$, then there are exactly 2 possibilities
\begin{enumerate}
	\item[1)] $\alpha_x=0$ and then $U^*_g (f)=0$ at the point $x$ since $c\neq 0$.
	\item[2)] Otherwise $\alpha_x\neq0$ and  then 
	$$-\frac{c}{n}\Big((n-1) \Delta f(x) -f(x) \Scal^g(x)\Big) \alpha_x = \d ^\nabla \Ric^g (\alpha,\alpha,\cdot)_x=0 \ , $$
which entails $(n-1) \Delta f(x) -f(x) \Scal^g(x)=0$ and so $U^*_g (f)=0$ at the point $x$.
\end{enumerate}
We can apply Theorem~\eq{KID} to conclude since $U^*_g (f)=0$ on the whole $M$.\qed

\begin{Rq}
In \cite{L} Section E1, Lafontaine gives a family of non trivial solutions of the equation $\hess^g f + \frac{\Delta f}{n}g - f\Ric^g_0=0 $. These metrics are again warped product metrics  $\d t^2 + h^2(t) g_0$ on $\Bbb S^1 \times Y$, where $(Y,g_0)$ is Einstein with positive scalar curvature, and $h$ a periodic function. Theorem~\eq{KID2} shows that the non trivial solutions of $(\Sigma_2)$ are the metrics in the class exhibited by Lafontaine in \cite{L}, that have positive constant scalar curvature. 
\end{Rq}

The case of surfaces is easy to treat since the computations of Lemmas~\eq{lem1}--\eq{lem2} and Theorem~\eq{main} are valid in dimension 2.

\begin{coro}
	Let  $(M^2, g)$ be a compact and connected Riemmanian surface which has a non trivial solution $(f,\alpha)$ of
$$ \left\{
\begin{array}{l}
	 \nabla\alpha + cfg= 0 \\
	 U_g^* (f)=\frac{1}{2}\Big( \Delta f -f \Scal^g\Big) g
\end{array}\right. \ ,$$
for a certain constant $c \neq 0$. Then  $(M^2,g)$ is isometric to the standard sphere $\mathbb S^2$.
\end{coro}
\proof The computations in the proof of Theorem~\eq{KID2} are still valid in dimension 2 so that $U^*_g (f)=0$ on the whole $M$. We can conclude using the surface version of the classification result on the system $(\Sigma_1)$.\qed

In the KID context, Theorem~\eq{KID2} gives rise to the natural issue\\
\question Let $(M^n,g)$ be a compact and connected manifold. Does the existence of non trivial solutions (i.e. $f\not\equiv0$) of 
$$
\left\{
\begin{array}{l}
	 \mycal{L}_\alpha g + 2cfg= 0 \\
	 U_g^* (f)=\frac{1}{n}\Big((n-1) \Delta f -f \Scal^g\Big) g
\end{array}\right. ,$$
(with $c$ a non zero constant) characterize the spatial Schwarzschild--de Sitter metrics for $n\ge 3$ (respectively the spatial de Sitter metrics for $n=2$)? \\


\subsection{The system $(\Sigma_3)$} The next main result deals with the restriction of the constraints application to  $\mycal C \times \Gamma(S^2 T^*M)$ where  $\mycal C= \left\{g\in\mycal M/ \d\Scal^g=0 \textrm{ and }   \Vol_{g}=\Vol(\Sn)\right\} $. In this context, $U^*_g$ has to be equal to the traceless Ricci curvature, so that we need to consider the problem of finding some $(f,\alpha)$ such that
$$\left\{
\begin{array}{l}
	 \mycal{L}_\alpha g + 2cfg= 0 \\
	 U_g^* (f)=\Ric^g_0
\end{array}\right. \ .$$
We could assume without loss of generality that $\Scal^g$ is a constant, but it is in fact not necessary in virtue of the following general result (which is the analogous version of Theorem~\eq{cmc} for the restriction to $\mycal C \times \Gamma(S^2 T^*M)$).

\begin{thm}\label{cmc2}
	Let  $(M^n, g)$ be a compact and connected manifold $(n\ge2)$ with $C^3$ metric, which has a solution $(f,\alpha), f\not\equiv0$, of
$$ \left\{
\begin{array}{l}
	 \mycal{L}_\alpha g+ 2cfg= 0 \\
	 U_g^* (f)=\Ric^g_0
\end{array}\right. \ ,$$
for a certain $c\in C^{\infty}(M), c\not\equiv 0$. Then the scalar curvature $\Scal^g$ is a positive constant and $c$ is a non zero constant.
\end{thm}
\proof We first prove that $\Scal^g$ is a positive constant. By computing the divergence of $ U_g^* (f)=\Ric^g_0$ we get
$\frac{1}{2}f \d \Scal^g=\delta U_g^* (f)=\delta\Ric^g_0=\left(\frac{1}{n}-\frac{1}{2}\right)\d\Scal^g $, which leads to the identity $\left(f-\frac{2-n}{n}\right)\d\Scal^g=0$. Let us define $F$ a closed subset as $F:=f^{-1}\left(\frac{2-n}{n}\right)$, and prove by contradiction that its interior $\Int(F)$ is empty. If $\Int(F)\neq\emptyset$ then, on $\Int(F)$  we have $U_g^* (f)=-\left(\frac{2-n}{n}\right)\Ric^g=\Ric^g_0$ whose trace gives $\Scal^g=0$. But on its complement $M\smallsetminus \Int(F)$, we have $\d\Scal^g=0$ and so $\Scal^g=0$ on the whole $M$ by a continuation argument. The first variation Formula~(\eq{t5}) of the Ricci curvature applies since $\mycal{L}_\alpha g=- 2cfg$:
$$ \nabla_\alpha \Ric^g (X,Y)= -(n-2)\hess^g (cf)(X,Y) - 2cf \Ric^g (X,Y) + (\Delta (cf))g(X,Y) \ ,$$
whose trace is zero since $\Scal^g\equiv0$ is a constant. We derive $\Delta (cf)=\frac{\Scal^g}{n-1}cf=0$, which implies that the function $cf$ is a constant. This is in contradiction with the fact that $\alpha$ is non isometric and $f\not\equiv0, c\not\equiv0$. Therefore, $\Int(F)=\emptyset$ and so $M\smallsetminus F$ is an open and dense subset of $M$ where $\d\Scal^g=0$. We finally obtain that $\Scal^g$ is a nonnegative constant (since $\Delta (cf)=\frac{\Scal^g}{n-1}cf$) which cannot be zero ($\alpha$ is non isometric).\\
The trace of $ U_g^* (f)=\Ric^g_0$ gives $\Delta f=\frac{\Scal^g}{n-1}f$ so the functions $f$ and $cf$ belong to the same eigenspace of positive eigenvalue $\frac{\Scal^g}{n-1}$. Unfortunately, $f^{-1}(0)$ the nodal set of $f$ has possibly a singular (or degenerate) part since the argument of Theorem~\eq{Bourg} does not work anymore. We need to use a crucial result on nodal sets owed to Bär, namely Corollary~2 in \cite{Baer}. It states that the nodal set of $f$ is the disjoint union $f^{-1}(0)=N_{reg}\cup N_{sing}$ with $N_{reg}:=\left\{x\in f^{-1}(0)/ \d_x f\neq 0\right\}$ and $N_{sing}:=\left\{x\in f^{-1}(0)/ \d_x f = 0\right\}$ that have the following properties:
\begin{enumerate}
	\item $N_{reg}$ is composed with smooth embedded hypersurfaces (by the implicit function theorem),
	\item $N_{sing}$ is a countably $(n-2)$--rectifiable set and thus has Hausdorff dimension $n-2$ at most.
\end{enumerate}
We notice that $N_{reg}$ cannot be empty. Indeed (by contradiction) if it was empty then there should be a unique nodal domain $M\smallsetminus N_{sing}$ (since $N_{sing}$ has Hausdorff dimension $n-2$ at most, $M\smallsetminus N_{sing}$ is connected). But $f$ has a vanishing integral on $M$ and has a constant sign on the dense open set $M\smallsetminus N_{sing}$, thereby $f\equiv 0$ which is not possible. Now consider $\Omega$ any connected component of the open set $M\smallsetminus N_{reg}$. By construction the boundary $\partial\Omega$ is smooth and $f$ is a solution of the usual Dirichlet problem on $\Omega$ since by construction we have 
$$\left\{
\begin{array}{cl}
	\Delta f=\frac{\Scal^g}{n-1}f& \textrm{on } \Omega  \\
	f=0 & \textrm{on } \partial\Omega
\end{array}\right. \ ,$$
with positive eigenvalue $\frac{\Scal^g}{n-1}>0$. But $f$ has a constant sign on $\Omega$ ($f$ can possibly vanish but only on a set of Hausdorff dimension less or equal to $n-2$) and so $\frac{\Scal^g}{n-1}>0$ has to be the first Dirichlet eigenvalue of the domain $(\Omega,g)$ which is always simple. Now, the function $cf$ also satisfies
$$\left\{
\begin{array}{cl}
	\Delta (cf)=\frac{\Scal^g}{n-1}cf& \textrm{on } \Omega  \\
	cf=0 & \textrm{on } \partial\Omega
\end{array}\right. \ ,$$
and consequently, there exists a non zero (since both $f$ and $cf$ cannot vanish on the whole $\Omega$) constant $\lambda$ such that $cf=\lambda f$. We have then proved that $c$ has to be a non zero constant on (a dense open subset of) $\Omega$. This  argument is true on any connected component of the open set $M\smallsetminus N_{reg}$ whose union is a dense open set in $M$, which means that $c$ has to be a non zero constant on the whole $M$ by a continuation argument.\qed

As a consequence we can completely answer the question of characterizing the manifolds that have non trivial solutions of 
$$ \left\{
\begin{array}{l}
	 \mycal{L}_\alpha g+ 2cfg= 0 \\
	 U_g^* (f)=\Ric^g_0
\end{array}\right. \ ,$$
in dimension 2  and 3.
\begin{coro}
Let  $(M^n, g)$ be a compact and connected Riemannian manifold of dimension $n=2$  or $3$, which has a solution $(f,\alpha), f\not\equiv0$, of
$$ \left\{
\begin{array}{l}
	 \mycal{L}_\alpha g+ 2cfg= 0 \\
	 U_g^* (f)=\Ric^g_0
\end{array}\right. \ ,$$
for a certain $c\in C^{\infty}(M), c\not\equiv 0$. Then $c$ is a non zero constant and $(M^n,g)$ is isometric to the standard sphere $\mathbb S^n$. From the KID point of view, the slice $(M^n, g,  cg)$ is a spatial de Sitter metric.
\end{coro}
\proof From Theorem~\eq{cmc2} we know that $c$ is a non zero constant and $\Scal^g$ a positive constant.\\
For the case of Riemannian surfaces i.e. $n=2$, Gauss--Bonnet formula claims that $M^2$ is homeomorphic to $\mathbb S^2$ and Theorem 3.83 in \cite{GHL} says that $M^2$ is in fact isometric to a standard 2--dimensional sphere (since the scalar curvature and the sectional curvature are equivalent for surfaces). The volume  normalization says that it is exactly $\mathbb S^2$.\\
When the dimension is $n=3$, then the theorem of Bessières--Lafontaine--Rozoy applies and so $M^3$ is isometric to $\mathbb S^3$.\qed

For dimension greater or equal to 4, the problem is quite harder to solve. Here again, we study this system assuming the closed character of the Killing conformal form $\alpha$.

\begin{thm}\label{BLR}
Let $(M^n, g)$ be a compact and connected manifold $(n\ge 2)$ with $C^3$ metric. Suppose there exists a solution $(f,\alpha), f\not\equiv0$, of
$$(\Sigma_3)\quad \left\{
\begin{array}{l}
	 \nabla\alpha+ cfg= 0 \\
	 U_g^* (f)=\Ric^g_0
\end{array}\right. \ ,$$
$c\in C^{\infty}(M), c\not\equiv 0$. Then c is a non--zero constant and $(M^n,g)$ is isometric to the standard sphere $\Sn$. From the KID point of view, the slice $(M^n, g,cg)$ is a spatial de Sitter metric.
\end{thm}

\proof Thanks to Theorem~\eq{cmc2} we know that $c$ is a non zero constant and that $\Scal^g$ is a positive constant. Thus we can use Theorem~\eq{main} so as to get
\begin{eqnarray*}
	f\d ^\nabla U_g^* (f) &=& \d f \wedge U_g^*(f) -f^2 \d ^\nabla\Ric^g \\
	&=&\d f \wedge \Ric^g_0 -f^2 \d ^\nabla\Ric^g 
\end{eqnarray*}
whereas a straight computation leads to $f\d ^\nabla U_g^* (f)= f \d ^\nabla \Ric^g_0= f \d ^\nabla \Ric^g $ since the scalar curvature is a constant. We identify both terms to obtain the crucial formula
$$ f(f+1) \d ^\nabla\Ric^g = \d f \wedge\Ric^g_0 \ .$$
Writing down Equation~(\eq{t6}) with $\psi =cf$ and taking our new information into account 
\begin{eqnarray*}
	-cf(f+1)\Ric_0^g &=& \d f(\alpha) \Ric^g_0 - \d f \otimes \Ric^g_0 (\alpha) \\
	&=& \d f(\alpha) \Ric^g_0 -c(n-1)\d f\otimes\d f + \frac{\Scal^g}{n}\d f\otimes \alpha  \ ,
\end{eqnarray*}
the trace of which gives $\d f(\alpha)\Scal^g=c n(n-1)\left|\nabla f\right|^2$. Then we decompose $M$ in the disjoint union $M=\Omega\coprod \mathcal{C}$, where we define the open set $\Omega=\left\{x\in M/ \d_x f\neq 0\right\}$ and $\mathcal{C}=\Omega^c$ is the set of the critical points of $f$. The closed set $\mathcal{C}$ possibly has a non empty interior, and then we write  $M=\Omega\coprod \partial\mathcal{C}\coprod \overset{\circ}{\mathcal{C}}$. By construction the open set $\Omega\coprod\overset{\circ}{\mathcal{C}}$ is dense in $M$. For any $x\in \Omega\coprod\overset{\circ}{\mathcal{C}}$ we have the 2 possibilities:
\begin{enumerate}
	\item  $x\in \overset{\circ}{\mathcal{C}}$. We are going to see that $f=0$ on the open set $\overset{\circ}{\mathcal{C}}$. Indeed, by definition $\d f=0$ on $\overset{\circ}{\mathcal{C}}$ which means that $f$ is locally constant and we have $\Delta f=0=\frac{\Scal^g}{n-1}f$, with $\Scal^g>0$ and thereby $f=0$ (we have also $\hess^g f=0$). Thus, the equation $U^*_g (f)=\Ric^g_0$ along $\overset{\circ}{\mathcal{C}}$ gives $\Ric^g_0=0$ (in particular $(\Ric^g_0)_x=0$).\\
	\item $x\in \Omega$ and in that case we denote by $\Omega_x$ the connected component of $\Omega$ that contains $x$. Thanks to Equation~(\eq{t2}), there exists an open ball $B(x,r) \subset\Omega_x$ $(r>0)$ where we can write ${\alpha=a(f)\d f}$ for a certain function $a\in C^{\infty}(\R)$. This implies $a(f)\Scal^g =cn(n-1)$ and so $a(f)$ is a non zero constant on $B(x,r)$. Since this argument is valid on a neighborhood of each point of $\Omega_x$, we get that $a(f)$ is a non zero constant on $\Omega_x$. We now use the first equation of $(\Sigma_3)$ (all the computations are carried out on the open set $\Omega_x$)
$$ \nabla\alpha=-cfg= \nabla\left(\frac{cn(n-1)}{\Scal^g}\d f\right)= \frac{cn(n-1)}{\Scal^g}\hess^g f \ , $$
that we plug in the second equation of $(\Sigma_3)$
\begin{eqnarray*}
	\Ric^g_0&=& \hess^g f -f\Ric^g + \frac{\Scal^g}{n-1}g\\
	&=& -f\Ric^g + \left(\frac{\Scal^g}{n-1}- \frac{\Scal^g}{n(n-1)}\right)g\\
	&=& -f\Ric^g_0 \ ,
\end{eqnarray*}
that is to say $(1+f)\Ric^g_0=0$ on $\Omega_x$. Particularly,
\begin{eqnarray*}
	0=-f(1+f)\Ric^g_0&=& \frac{cn(n-1)}{\Scal^g}\left|\nabla f\right|^2 \Ric^g_0 -c(n-1)\d f\otimes\d f + \frac{a(f)\Scal^g}{n}\d f\otimes \d f \\
	&=&\underbrace{ \frac{cn(n-1)}{\Scal^g}\left|\nabla f\right|^2}_{\neq 0} \Ric^g_0 \ ,
\end{eqnarray*}
that allows us to conclude $ \Ric^g_0=0$ on $\Omega_x$ ($\Rightarrow(\Ric^g_0)_x=0$).
\end{enumerate}
We have proved that $\Ric^g_0\in C^1$ (since $g\in C^3$) is identically zero on the dense open set $\Omega\coprod\overset{\circ}{\mathcal{C}}$, which induces $\Ric^g_0 =0$, namely $g$ is Einstein. But, it is well known thanks to Obata \cite{O}, that $(M^n,g)$ is isometric to the standard sphere $\Sn$.\qed

Theorem~\eq{BLR} must be related to the result in dimension 3 of Bessières--Lafontaine--Rozoy \cite{BLR} evocated in the introduction. Their proof which is specific to the dimension 3, consists on proving that any metric $g$ of constant scalar curvature having a solution $f$ of $U^*_g(f)=\Ric_0^g$, is conformally flat. They can conclude, by using a result of Lafontaine on conformally flat manifolds (\cite{L} Section E2), that $g$ is the standard metric of the sphere $\Bbb S^3$. It is not obvious that their result could be extended to higher dimensions, and this is the reason why Theorem~\eq{BLR} can be considered as a generalization of their theorem for any dimension $n\ge3$.\\
It is also of great interest that the existence of non trivial solutions of $(\Sigma_3)$ characterizes the standard spherical slices of de Sitter space--time. We then address the following\\
\question Let $(M^n,g)$ be a compact and connected manifold, $n\ge 4$. Does the existence of non trivial solutions (i.e. $f\not\equiv 0$) of 
$$\left\{
\begin{array}{l}
	 \mycal{L}_\alpha g + 2cfg= 0 \\
	 U_g^* (f)=\Ric^g_0
\end{array}\right. \ ,$$
characterize the geometry of the standard sphere $\Sn$?


\section{Another Equation}\label{4}

A slightly different point of view is to consider $U^*_g$ as a linear order 2 differential operator on functions. Naturally, one could think of modifying this operator thanks to a potential i.e. considering an equation of the kind $U^*_g(f)=f \tau$ where $\tau$ is a symmetric 2--tensor on $M$ ($\tau$ is the potential). Of course, the motivation of such an equation vanishes (the choice of $\tau$ is not clearly suggested by a geometric formulation) except the origin of the operator $U^*_g$ itself via the scalar curvature application. However, the situation that is the closest to the equation \textit{à la Bessières--Lafontaine--Rozoy} is to take $\tau=\Ric_0^g$, leading to the modified equation $U^*_g(f)=f \Ric^g_0$. In the continuation of the previous results of this article we consider the problem of finding some couple $(f,\alpha)$ such that
$$(\Sigma'_4)\quad \left\{
\begin{array}{l}
	 \mycal L_\alpha g+ 2cfg= 0 \\
	 U_g^* (f)= f \Ric^g_0
\end{array}\right. \ ,  $$
where $c\in C^\infty (M)$. The problem is that too many manifolds do have no trivial solutions of $(\Sigma'_4)$:
\begin{enumerate}
	\item [1)] Examples with $c=0$:
\begin{enumerate}
	\item A manifold $(M^n,g)$ that has a Killing form $\alpha$. Then the couple ${(f=0,\alpha)}$ is a solution of $(\Sigma'_4)$. The scalar curvature of $g$ is not necessarily a constant and note that this class of manifolds is quite (in fact too) large.
	\item A Ricci flat manifold $(M^n,g)$ that has a Killing form $\alpha$. Then any couple $(f=const.\neq 0,\alpha)$ is a solution of $(\Sigma'_4)$. The scalar curvature of $g$ is zero by construction and this class of manifolds is also quite large, since it contains the flat tori.
	\item The standard sphere $\Sn$. Then any couple $(f,\alpha)$ where $f\in \Ker U^*_g$ and $\alpha$ a Killing form, is a solution of $(\Sigma'_4)$.
\end{enumerate}
	\item [2)] Examples with $c\not \equiv 0 $:
\begin{enumerate}
	\item A manifold $(M^n,g)$ that has a Killing form $\alpha$. Then the couple $(f=0,\alpha)$ is a solution of $(\Sigma'_4)$ for any function $c\in C^\infty (M)$ (the non zero constant functions $c$ are obviously admissible). Once again the scalar curvature is not necessarily a constant and this class of manifolds is very large.
	\item The standard sphere $\Sn$. Then any couple $(f,c \d f)$ where $f\in \Ker U^*_g$ et and $c$ a nonzero constant, is a solution of $(\Sigma'_4)$.
\end{enumerate}
\end{enumerate}

These examples explain why we choose to restrict our problem to the case where $c$ is a nonzero constant and the Killing conformal form $\alpha$ is closed. Under these conditions, Theorem~\eq{BLRm} claims that the standard sphere $\Sn$ is the only manifold carrying non trivial solutions of our new system that we call $(\Sigma_4)$. For sake of simplicity, we normalize the non zero constant $c=1$ till the end of this section.

\begin{thm}\label{BLRm}
Let $(M^n, g)$ be a compact and connected manifold with $C^3$ metric $(n\ge2)$, which has a non trivial solution $(f,\alpha)$ of
	$$(\Sigma_4)\quad \left\{
\begin{array}{l}
	 \nabla\alpha+ fg= 0 \\
	 U_g^* (f)= f \Ric^g_0
\end{array}\right. $$	
Then $(M^n,g)$ is isometric to the standard sphere $\Sn$. 
\end{thm}

\proof The aim is to prove that $g$ is Einstein. We still have  $\Delta f= \frac{\Scal^g}{n-1}f$ by tracing $U_g^* (f)= f \Ric^g_0$. On the one hand, Equation~(\eq{eq}) of Lemma~\eq{lem2} entails
\begin{eqnarray*}
	f\d ^\nabla U_g^* (f) &=& \d f \wedge U_g^*(f) -f^2 \d ^\nabla\Ric^g +\Big( f\d (\Delta f) - (\Delta f) \d f\Big)\wedge g \\
	&=&f\d f\wedge \Ric^g_0 -f^2 \d ^\nabla\Ric^g +\left( f\d \left(\frac{\Scal^g}{n-1}f\right) -  \frac{\Scal^g}{n-1}f\d f\right)\wedge g\\
	&=&f\d f\wedge \Ric^g_0 -f^2 \d ^\nabla\Ric^g  + \frac{f^2}{n-1}\d\Scal^g \wedge g \ ,
\end{eqnarray*}
and on the other hand, we straightly compute
$$f\d ^\nabla U_g^* (f) =f\d f \wedge\Ric^g_0 + f^2 \d ^\nabla\Ric^g  - \frac{f^2}{n}\d\Scal^g \wedge g   \ .$$
We identify both expressions and get
$$ 2 f^2 \d ^\nabla\Ric^g = f^2 \left(\frac{2n-1}{n(n-1)}\right)\d\Scal^g \wedge g \ .  $$
Thanks to Equation~(\eq{t7}) it comes $f^2 \d\Scal^g \wedge \alpha =0$, but since the open set $f^{-1}(\R^*)$ is dense  (the argument of Bourguignon in \cite{B} works in our situation, see also Theorem~\eq{Bourg}), we  have $ \d\Scal^g \wedge \alpha =0$ on the whole $M$, because this 2--form is $C^0$ ($g\in C^3$). The trace of Equation~(\eq{t5}) gives
$$\d\Scal^g (\alpha)= 2f\Scal^g - 2(n-1)\Delta f =0 \ .$$
Consider Equation~(\eq{t6})
$$-2f^3\Ric^g_0= 2 f^2\d ^\nabla\Ric^g (\alpha,X,Y)= -f^2 \left(\frac{2n-1}{n(n-1)}\right)\d \Scal^g \otimes \alpha (X,Y) \ ,$$
that we shortly write 
$$ 2f^3\Ric^g_0=f^2 \left(\frac{2n-1}{n(n-1)}\right)\d \Scal^g \otimes \alpha \ .$$
If one takes  $x\in f^{-1}(\R^*)$, then there are 2 possibilities:
\begin{enumerate}
	\item[1)] The linear form $\alpha_x\equiv0$ and thereby $(\Ric^g_0)_x=0$.
	\item[2)] Otherwise the linear form $\alpha_x\not\equiv 0$. But we know  $ (\d\Scal^g \wedge \alpha)_x =0$ and also ${\d_x\Scal^g(\alpha_x)=0}$,  so that $\d_x\Scal^g=0$ and consequently $(\Ric^g_0)_x=0$.
\end{enumerate}
We obtain $\Ric^g_0=0$ on $M$ (by a continuation argument), namely $g$ is Einstein. But it is well known thanks to Obata \cite{O}, that $(M^n,g)$ is isometric to the standard sphere $\Sn$.\qed 

Theorem~\eq{BLRm} appears as a good generalization of the Obata equation since the standard sphere is the unique manifold that has a non trivial solutions to $(\Sigma_4)$. This result is very general since we do not need to make any assumption on the scalar curvature and since it is valid in any dimension.\\


\end{document}